\documentclass[a4paper,fleqn,11pt]{amsart}

\usepackage{mathrsfs}
\usepackage{natbib}
\usepackage{hyperref}
\hypersetup{
	pdftitle={L\'evy Processes in the Heavy-Traffic Regime},
	pdfsubject={Probability Theory},
	pdfauthor={Kamil Kosi\'nski},
	pdfdisplaydoctitle=true,
}

\newtheorem{theorem}{Theorem}
\theoremstyle{plain}
\newtheorem{proposition}{Proposition}
\theoremstyle{remark}
\newtheorem{remark}{Remark}

\newcommand{\pp}{\mathbb{P}}
\newcommand{\ee}{\mathbb{E}}
\newcommand{\rr}{\mathbb{R}}

\newcommand{\ML}{\mathscr{ML}}
\newcommand{\RV}{\mathscr{RV}}
\newcommand{\LA}{\mathscr{L}^{(\alpha)}}

\newcommand{\dto}{\stackrel{d}{\to}}
\newcommand{\de}{\stackrel{d}{=}}

\newcommand{\toi}{\to\infty}
\newcommand{\as}[1]{\quad\text{as}\quad#1\toi}
\newcommand{\bx}[1]{\bar X^{(#1)}}
\newcommand{\bs}[1]{\bar S^{(#1)}}

\newcommand{\lr}{\mathscr R}

\begin{document}
\bibliographystyle{plainnat}
\setcitestyle{numbers}

\title[Heavy-traffic]{Convergence of the all-time supremum of a L\'evy process
in the heavy-traffic regime}

\author{K.M.\ Kosi\'nski}
\address{E{\sc urandom}, Eindhoven University of Technology, P.O. Box 513, 5600 MB  Eindhoven, the Netherlands \\
Korteweg-de Vries Institute for Mathematics, University of Amsterdam, P.O. Box 94248, 1090 GE Amsterdam, the Netherlands}
\email{Kosinski@eurandom.tue.nl}
\thanks{The first author was supported by NWO grant 613.000.701.}

\author{O.J.\ Boxma}
\address{E{\sc urandom} and Department of Mathematics and Computer Science \\
Eindhoven University of Technology, P.O. Box 513, 5600 MB  Eindhoven, the Netherlands}
\email{Boxma@win.tue.nl}

\author{B.\ Zwart}
\address{CWI, P.O. Box 94079, 1090 GB Amsterdam, the Netherlands}
\email{Bert.Zwart@cwi.nl}

\date{\today}

\begin{abstract}
In this paper we derive a technique of obtaining limit theorems for suprema of L\'evy processes
from their random walk counterparts. For each $a>0$, let $\{Y^{(a)}_n:n\ge 1\}$
be a sequence of independent and identically distributed random variables
and $\{X^{(a)}_t:t\ge 0\}$ be a L\'evy process such that $X_1^{(a)}\de Y_1^{(a)}$, $\ee X_1^{(a)}<0$
and $\ee X_1^{(a)}\uparrow0$ as $a\downarrow0$. Let $S^{(a)}_n=\sum_{k=1}^n Y^{(a)}_k$.
Then, under some mild assumptions,
$\Delta(a)\max_{n\ge 0} S_n^{(a)}\dto \lr\iff\Delta(a)\sup_{t\ge 0} X^{(a)}_t\dto \lr$, 
for some random variable $\lr$ and some function $\Delta(\cdot)$. We utilize this result
to present a number of limit theorems for suprema of L\'evy processes in the heavy-traffic regime.
\end{abstract}

\maketitle

\section{Introduction}
For each $a>0$, let $X^{(a)}\equiv\{X^{(a)}_t:t\ge 0\}$ be a L\'evy process
such that $\mu^{(a)}:=\ee X^{(a)}_1<0$. Along with the L\'evy process $X^{(a)}$
define $\bx a=\sup_{t\ge 0} X^{(a)}_t$. Since $\mu^{(a)}<0$ assures that
$X^{(a)}$ drifts to $-\infty$, the {\it all-time} supremum $\bx a$ is a proper
random variable for each $a>0$. Now if $\mu^{(a)}\uparrow 0$ as $a\downarrow 0$, then 
$\bx a\toi$. From this fact
a natural question arises: How fast does $\bx a$ grow as $a\downarrow 0$? 

The main purpose of this paper is to answer the above question by considering
the discrete approximation of a L\'evy process by a random walk. 
For each $a>0$, let $\{Y_n^{(a)}:n\ge1\}$ be a sequence of independent and identically distributed (i.i.d.) random variables. Put 
$\bs a=\max_{n\ge 0} S_n^{(a)}$, 
where $S_n^{(a)}=\sum_{k=1}^n Y_k^{(a)}$, $S^{(a)}_0=0$.
We shall show that if $Y_1^{(a)}$ has the same distribution as $X_1^{(a)}$, then the limiting
distribution of $\bx a$ can be derived from the limiting distribution of $\bs a$.
In doing so we shall utilize a bound by \citet{Willekens87}.
Loosely speaking, this bound
allows to derive certain properties of L\'evy processes via their corresponding
random walk approximations (see also \citet{Doney04}). The advantage of this approach is that
the problem on how fast does $\bs a$ grow as $a\downarrow 0$ has been treated 
extensively and various methods have been developed.

One major reason why the behaviour of $\bs a$ has been studied is that it is well known 
that the stationary distribution of the waiting time of a customer
in a single-server first-come-first-served $GI/GI/1$ queue coincides with the distribution of the 
maximum of a corresponding random walk. 
The condition on the mean of the random walk becoming small ($a\downarrow 0$)
means in the context of a queue that the traffic load tends to 1. Thus, the problem 
under consideration (in the random walk setting)
may be seen as the investigation of the growth rate of the stationary waiting-time distribution in a $GI/GI/1$
queue. This is one of the most important problems in 
queueing theory that is referred to as the heavy-traffic approximation problem. 
The question was first posed by Kingman (see \cite{Kingman65} for an extensive discussion 
on the early results). It has been solved in various settings by, e.g.,
\citet{Prokhorov63}, \citet{Boxma99}, \citet{Resnick00}, \citet{Szczotka04} and many others.

Surprisingly, there are no results in
the literature on the heavy-traffic limit theorems for L\'evy-driven ({\it fluid}) queues.
Our approach however allows to translate each single result in the random walk setting
to its analogue in the L\'evy setting, thereby providing a range of fluid heavy-traffic
limit theorems. Our main result, \autoref{thm:main},
states: under some mild conditions,  
$\Delta(a)\bs a \dto \lr$ if and only if $\Delta(a)\bx a \dto \lr$, 
for some random variable $\lr$ and some function $\Delta(\cdot)$.

The remainder of the paper is organized as follows. 
In \autoref{PRE} we fix notation and give some necessary preliminaries.
\autoref{MR} contains the main result of this paper, \autoref{thm:main}, and its proof.
Instances of this theorem applied to the results by \citet{Boxma99}, \citet{Shneer09}
and \citet{Szczotka03} (see also \citet{Czystolowski10}) are presented in \autoref{SI} and conclude the paper.

\section{Preliminaries and notation}
\label{PRE}
All the stochastic objects are assumed to be defined on the probability space
$(\Omega,\mathbb F,\pp)$ endowed with a standard filtration $\mathcal F=\{\mathcal F_t:t\ge0\}$,
that is $\mathcal F$ is an increasing, right-continuous family of complete sub-$\sigma$-fields
of $\mathbb F$. Let us begin by fixing the notation for L\'evy processes.

A real-valued stochastic process $X\equiv\{X_t:t\ge0\}$, with $X_0=0$, 
is said to be a L\'evy process with respect to the filtration $\mathcal F$ if 
it is adapted to $\mathcal F$, $X_s - X_t$ is independent of $\mathcal F_t$ and distributed
as $X_{s-t}$ for any $0\le t<s$. Moreover, we assume that the sample paths of $X$
are {\it c\`adl\`ag} (right-continuous with left limits), so that $X$ is strong Markov.

Let $\psi$ be the L\'evy characteristic exponent of $X$ so that $\ee e^{iuX_t}=e^{-t \psi(u)}$, for all $u\in\rr$. 
In this case, for some $\sigma>0$ and $\delta\in\rr$, $\psi$ has the form
\[
\psi(u)=i\delta u+\frac{1}{2}\sigma^2u^2+\int_{|x|<1}\left(1-e^{iu x}+iu x\right)\nu(dx)
+\int_{|x|\ge 1}\left(1-e^{iu x}\right)\nu(dx),
\]
where $\nu$ is the L\'evy measure (on $\rr\setminus\{0\}$) satisfying $\int_{\rr}(1\wedge x^2)\nu(dx)<\infty$.
We only consider {\it nondeterministic} processes $X$, which is synonymous with $\sigma^2+\nu(\rr\setminus\{0\})>0$. 
We say that: $X$ is {\it centered} if $\ee X_t=0$ for all $t$; 
{\it spectrally positive} if $\nu(-\infty,0)=0$; {\it spectrally negative} if $\nu(0,\infty)=0$.
If $X_1$ has a stable distribution with index $\alpha\in(0,2]$ then we say that $X$ is 
an $\alpha$-{\it stable} L\'evy process and denote it by $\LA$.
For more background on L\'evy processes we refer the reader to \cite{Bertoin90} and references therein.

For an interval $B\subseteq[0,\infty)$, we denote by $D(B)$ the space of real-valued c\`adl\`ag functions on $B$
equipped with the usual Skorokhod $J_1$ topology; see, e.g., \cite[Chapter VI]{Jacod87}.

In the sequel we will encounter the Mittag--Leffler distribution; see,
e.g., \cite[p.\ 329]{Bingham87}.
A positive random variable $\mathscr{ML}_\alpha$ is said to have a Mittag--Leffler distribution with parameter $\alpha\in(0,1]$
if its Laplace--Stieltjes transform (LST) is given by
\[
\ee \exp(-s \mathscr{ML}_\alpha) =\frac{1}{1+s^{\alpha}}.
\]
Observe that $\mathscr{ML}_1$ has the 1-exponential distribution.

We will also make use of some standard notation. For two functions $f$, $g$ we shall
write $f(x)\sim g(x)$ as $x\to x_0\in[0,\infty]$ to mean $\lim_{x\to x_0}f(x)/g(x)=1$. 
The class of regularly varying functions with index $\alpha$ shall be denoted
by $\mathscr{RV}_\alpha$.

In what follows we shall also write
\[
\bx a=\sup_{t\ge 0} X^{(a)}_t,\quad \bs a=\max_{n\ge 0} S^{(a)}_n,
\]
where, for each $a>0$, $\{X^{(a)}_t:t\ge 0\}$ is a L\'evy process and
$S^{(a)}_n=\sum_{k=1}^n Y_k^{(a)}$ is the $n$th partial sum of a sequence
of random variables $\{Y^{(a)}_n:n\ge1\}$ with $S^{(a)}_0=0$.
\section{Main theorem}
\label{MR}
\begin{theorem}
\label{thm:main}
For each $a>0$, let $\{Y^{(a)}_n:n\ge 1\}$ be a sequence of i.i.d. random variables and
$\{X^{(a)}_t:t\ge 0\}$ be a L\'evy process. 
Assume that for each $a$, $Y^{(a)}_1\de X^{(a)}_1$, $\mu^{(a)}=\ee Y_1^{(a)}<0$ and $\mu^{(a)}\uparrow 0$ as $a\downarrow 0$.
Then, for some random variable $\lr$, 
\begin{equation}
\label{eq:thm:main}
\Delta(a)\max_{n\ge 0} S_n^{(a)}\dto \lr\iff
\Delta(a)\sup_{t\ge 0} X^{(a)}_t\dto \lr,\quad\text{as}\quad a\downarrow 0,
\end{equation}
where $\Delta(\cdot)$ is a normalizing function such that $\Delta(a)X^{(a)}_1\to 0$ almost surely.
\end{theorem}
\begin{remark}
The assumptions of the above theorem are natural:
\begin{itemize}
\item $\mu^{(a)}\uparrow 0$ assures that, for each $a$, $\bx a$ and $\bs a$ are finite random variables.
Moreover,
$\bx a$ and $\bs a$ tend to infinity with $a$. Therefore the function $\Delta(\cdot)$ tends to zero and can be seen
as the speed of convergence in \eqref{eq:thm:main};
\item $\Delta(a)X^{(a)}_1\to 0$ is satisfied in all typical applications; 
for instance, when $X^{(a)}_t=X_t-at$ for a fixed L\'evy process $X$, see \autoref{subsec:SPP} and \autoref{RV};
or there exists a function $d(a)\toi$ such that $X^{(a)}_{d(a)t}\to X$ in $D[0,\infty)$, 
see \autoref{subsec:HTIP}.
\end{itemize}
The distribution of the random variable $\lr$ can be computed in several cases; we will get back to it in \autoref{rem:distr}.
\end{remark}
\begin{proof}[Proof of \autoref{thm:main}]
With $\bar R(x):=1-R(x)$, where $R$ is the distribution function of $\lr$, 
it is enough to show that, as $a\downarrow 0$,
\[
\pp\left(\Delta(a)\bs a>x\right)\to \bar R(x)\iff
\pp\left(\Delta(a)\bx a>x\right)\to \bar R(x),
\]
for any continuity point $x$ of $R$.

Assume that $\Delta(a)\bs a\dto\lr$, the converse implication follows
in the same manner.

Observe that $\bs a\de\max_{n\ge 0} X^{(a)}_n$.
Thus, we trivially have
\begin{equation}
\label{lower}
\pp\left(\Delta(a)\bx a > x \right)\ge
\pp\left(\Delta(a)\bs a > x \right).
\end{equation}
On the other hand, for any $x_0>0$,
\begin{align*}
\pp\left(\Delta(a)\bx a> x\right)
&\le 
\pp\left(\Delta(a)\max_{n\ge 0} X^{(a)}_n> x-x_0\right)
+
\pp\left(\Delta(a)\bx a> x,\,\Delta(a) \max_{n\ge 0} X^{(a)}_n\le x-x_0\right)\\
&=
\pp\left(\Delta(a)\bs a> x-x_0\right)
+
\pp\left(\Delta(a)\bx a> x,\,\Delta(a) \max_{n\ge 0} X^{(a)}_n\le x-x_0\right).
\end{align*}
We use an argument similar to the one in \cite{Maulik06, Willekens87}.
Define $\tau^{(a)}(x)=\inf\{t\ge 0:\Delta(a)X^{(a)}_t > x\}$, so that
$\tau^{(a)}$ is a stopping time \cite[Corollary 8]{Bertoin90}. Now
the second term on the right-hand side of the above inequality can be bounded from above by
\begin{align*}
\lefteqn{
\pp\left(\tau^{(a)}(x)<\infty,
\Delta(a)\left(\inf_{t\in[\tau^{(a)}(x),\tau^{(a)}(x)+1]} 
\left(X^{(a)}_t-X^{(a)}_{\tau^{(a)}(x)}\right)
\right)
\le-x_0\right)}
\\
&= 
\pp\left(\tau^{(a)}(x)<\infty\right)\pp\left(\Delta(a)\inf_{t\in[0,1]} X^{(a)}_t\le-x_0\right),
\end{align*}
where we used the strong Markov property in the last equality. Thus,
\begin{equation}
\label{upper}
\pp\left(\Delta(a)\bx a> x\right)\pp\left(\Delta(a)\inf_{t\in[0,1]} X^{(a)}_t >-x_0\right)
\le \pp\left(\Delta(a)\bs a> x-x_0\right).
\end{equation}
Now $\Delta(a)X_1^{(a)}\to0$ a.s. implies that the finite-dimensional distributions of
$\{\Delta(a)X_t^{(a)}:t\in[0,1]\}$ converge to zero a.s. and thus by \cite{Skorokhod57} the whole process
converges to 0 in $D[0,1]$. Applying the continuous mapping theorem with the infimum (over $[0,1]$) map
yields $\Delta(a)\inf_{t\in[0,1]} X^{(a)}_t\to 0$.
Thus, combining formulas \eqref{lower} and \eqref{upper} we get
\[
\bar R(x)\le \liminf_{a\to0}\pp\left(\Delta(a)\bx a> x\right)
\le\limsup_{a\to0}\pp\left(\Delta(a)\bx a> x\right)
\le\bar R(x-x_0).
\]
The thesis follows by letting $x_0\to 0$.
\end{proof}
\begin{remark}
In the following section,
we shall use the {\it if} part of \autoref{thm:main} to derive various limit theorems
for suprema of L\'evy processes. 
It is worth noting however that the {\it only if} part could
be used as well to derive limit theorems for suprema of random walks. A variation of this
approach has been undertaken in \citet{Szczotka03}, where first a heavy-traffic limit theorem is derived
in continuous time and then this theorem is used to claim an analogous behaviour in discrete time.
\end{remark}
\section{Special instances}
\label{SI}
\autoref{thm:main} provides a tool for translating limit theorems
for random walks to their analogues in the L\'evy setting. In this section
we shall focus our attention on some seminal results about
the convergence of the maxima of random walks and reformulate them
to the L\'evy case. We illustrate each special case that we consider
with a remark that explains an alternative way of obtaining the particular result
via a direct approach undertaken in the literature. These remarks, albeit short,
are rigorous enough to act as alternative proofs. Let us start
with the case in which the underlying processes are spectrally positive, which
is closely related to the queueing setting via the compound Poisson process.

\subsection{Spectrally positive processes}
\label{subsec:SPP}
For a sequence of zero mean, i.i.d.
random variables $\{Y,Y_n:n\ge1\}$, the question of how fast does
$\bs a=\max_{n\ge 0} (S_n-na)$ grow as $a\downarrow 0$
was first posed by \citet{Kingman61,Kingman62}. 
Kingman in his proof assumed exponential moments of $|Y|$ and
used Wiener--Hopf factorization to obtain the Laplace transform of $\bs a$.
\citet{Prokhorov63} generalized Kingman's result to the case when only the second moment of $Y$ is finite.
His approach was based on the functional Central Limit Theorem. These two approaches have become classical
and have both been used to prove various heavy-traffic results. The analytical approach of Kingman was used by
\citet{Boxma99} (see also \citet{Cohen02}) to study the limiting behaviour of $\bs a$ in the case of infinite variance.
They proved that if $\pp(Y>x)$ is regularly varying at infinity with a parameter $\alpha\in(1,2)$
(and under some additional assumptions), then there exists a function $\Delta(\cdot)$ such that
$\Delta(a)\bs a$ converges in law to a proper random variable. 

\cite[Theorem 5.1]{Boxma99} acts as the first application of our main result. 
For a L\'evy measure $\nu$ define
\[
r(s):=\int_0^\infty \left(e^{-sx}-1+sx\right)\,\nu(dx).
\]
For a L\'evy process $X\equiv\{X_t:t\ge 0\}$, let $F$ be the distribution function of $X_1$ and set $\bar F:=1-F$.
With $\alpha>0$,
\cite[Theorem 8.2.1]{Bingham87} asserts that $\bar F\in\RV_{-\alpha}$ if and only if $\nu(x,\infty)\in\RV_{-\alpha}$, where
$\nu$ is the L\'evy measure of $X$; moreover: $\bar F(x)\sim \nu(x,\infty)$, as $x\toi$. This, combined
with \cite[Theorem 5.1]{Boxma99} and \autoref{thm:main}, yields:
\begin{proposition}
\label{thm:Boxma}
Let $X$ be a spectrally positive L\'evy process such that $\nu(x,\infty)\in\mathscr{RV}_{-\alpha}$ for $\alpha\in(1,2)$. Set
$\rho(a)=\mu/a$, where $\mu=\ee X_1$, then
\[
\Delta(\rho(a))\sup_{t\ge0} (X_t- at)\dto \ML_{\alpha-1},\quad\text{as}\quad \rho(a)\uparrow 1,
\]
where $\Delta(x)=d(x)/\mu$ and $d(x)$ is such that
\begin{equation}
\label{eq:contraction}
r(d(x))\sim d(x)\frac{1-x}{x}\mu^{\alpha},\quad\text{as}\quad x\uparrow 1.
\end{equation}
\end{proposition}
See also \cite{Boxma99,Resnick00} for possible refinements of the assumption on regular variation in this special case.
\begin{remark}
\label{Pollaczek}
It is possible to prove \autoref{thm:Boxma} using a direct approach like the one in \cite{Boxma99}.
Let $X^{(a)}_t=X_t-at$, $a>\mu$, then the Pollaczek--Khinchine formula (see, e.g., \cite{Asmussen03} Chapter IX) yields
\[
\ee e^{-\lambda \bx a}
=\frac{\lambda\varphi_a'(0)}{\varphi_a(\lambda)}<\infty,\quad\text{for}\quad\lambda>0,
\]
where $\varphi_a(\lambda)=\log\ee\exp(-\lambda(X_1-a))$. 
Substituting $\varphi_a(\lambda)=\lambda\varphi_a'(0)+r(\lambda)$ yields
\[
\frac{\lambda\varphi_a'(0)}{\varphi_a(\lambda)}=\frac{1}{1+\frac{r(\lambda)}{\lambda\varphi_a'(0)}},
\]
where we assumed $\sigma=0$ for simplicity. Let $\lambda= s \Delta(\rho(a))$ 
with $\Delta(\cdot)$ as in \autoref{thm:Boxma}.
Using \cite[Theorem 8.1.6]{Bingham87} one infers that,
under the assumption $\nu(x,\infty)\in\RV_{-\alpha}$,
$r$ is a regularly varying function at 0 with index $\alpha$. 
We necessarily have $d(x)\downarrow 0$, as $x\uparrow 1$. Hence, as $\rho(a)\uparrow 1$,
\[
\frac{r(\lambda)}{\lambda\varphi_a'(0)}
\sim
\left(\frac{s}{\mu}\right)^{\alpha-1}
\frac{r(d(\rho(a)))}{d(\rho(a))(a-\mu)}
=
\frac{s^{\alpha-1}}{\mu^{\alpha}}
\frac{r(d(\rho(a)))}{d(\rho(a))}
\frac{\rho(a)}{1-\rho(a)}
\sim
s^{\alpha-1}.
\]
\end{remark}

\subsection{Asymptotically stable processes}
\label{RV}
\autoref{thm:Boxma} limits the class of L\'evy processes under consideration
to spectrally positive. Further improvements of the result from \cite{Boxma99} by
\citet{Furrer97} and \citet{Resnick00} assumed that the random walk belongs to the domain of
attraction of a spectrally positive stable law and relied on functional limit theorems. 
\citet{Shneer09} relaxed this assumption to allow the random walk to belong
to the domain of attraction of any stable law.
The main result from \cite{Shneer09} acts as the second instance of an application of \autoref{thm:main}.
\begin{proposition}
\label{thm:Seva}
Let $X$ be a centred L\'evy process such that the random variable 
$X_1$ belongs to the domain of attraction of a stable law $\LA_1$ with index $\alpha\in(1,2]$. That is,
there exists a sequence $\{d(n):n\ge 0\}$ such that
\begin{equation}
\label{ass:CLT}
\frac{X_n}{d(n)}\dto\LA_1,\as n.
\end{equation}
Then,
\[
\Delta(a)\sup_{t\ge 0} (X_t-at )\dto \sup_{t\ge0}(\LA_t-t),\quad\text{as}\quad a\downarrow 0,
\quad\text{where}\quad \Delta(a)=\frac{1}{d(n(a))}
\]
and $n(a)$ is such that
\begin{equation}
\label{eq:defna}
an(a)\sim d(n(a)),\quad\text{as}\quad a\downarrow 0.
\end{equation}
\end{proposition}

\begin{remark}
\label{rem:distr}
It is well known that the sequence $d(n)$ in \autoref{thm:Seva} is regularly varying with index
$1/\alpha$. Therefore, \autoref{thm:Seva}
implies that, with $X^{(a)}_t=X_t-at$, $\bx a$ 
grows as a regularly varying function (at zero) with index $-1/(\alpha-1)$.
If $\LA$ is spectrally negative, then the distribution of $\lr=\sup_{t\ge0}(\LA_t-t)$ is exponential;
see, e.g., \cite[Proposition 5]{Bingham75}. 
If $\LA$ is spectrally positive, then, as seen in  
\autoref{thm:Boxma}, the limiting random variable has a Mittag--Leffler distribution;
see, e.g., \cite[Theorem 4.2]{Kella91}. 
If $\LA$ is symmetric, then one can give the Laplace--Stieltjes transform of 
$\lr$,
see \cite[Theorem 8]{Szczotka03}. 
In the other cases the explicit form of the distribution might be infeasible to compute; 
however, one can easily find its tail 
asymptotics $\pp(\lr>x)\sim C x^{1-\alpha}$. 
For more details on the supremum distribution of a L\'evy process see \cite{Szczotka03}.
\end{remark}
\begin{remark}
\citet{Shneer09} showed that both classical approaches, i.e., 
via Wiener--Hopf factorization and via a functional central limit theorem,
can be applied to obtain their result. Moreover, the technical difficulties arising from these methods
can be overcome using a generalization of Kolmogorov's inequality based on a result by \citet{Pruitt81}.
A similar result is also available for L\'evy processes and can also be found in \cite{Pruitt81}. 
Let us introduce $V(x)=\int_{|y|\le x} y^2 \nu(dy)$, the truncated second moment of the L\'evy measure $\nu$.
Under the assumptions of \autoref{thm:Seva}, $V\in\RV_{2-\alpha}$. Moreover,
\cite[Section 3]{Pruitt81} asserts that
there exists a constant $C$ such that
\begin{equation}
\label{eq:Pruitt}
\pp\left(\sup_{s\le t} X_s\ge x\right)\le C\frac{t V(x)}{x^2}.
\end{equation}
Using the regular variation of $V$, \eqref{eq:defna} and \eqref{eq:Pruitt}, 
for any fixed $T>0$ there exist constants $C_1,C_2>0$ such that 
\begin{align}
\pp\left(\sup_{t\ge n(a) T} \left( X_t-at\right)\ge 0\right )
&\le 
\sum_{k=0}^\infty \pp\left(\sup_{t\le 2^{k+1} n(a) T} X_t\ge 2^k a n(a) T\right )\nonumber\\
\label{eq:estimate}
&\le C_1\frac{V(an(a)T)}{a^2n(a)T}\sum_{k=0}^\infty (2^k)^{1-\alpha}
\le C_2\frac{V(d(n(a)))}{c^2(n(a))}n(a) T^{1-\alpha}.
\end{align}
The sequence $d(n)$ can be defined as $d(n):=\inf\{t>0:V(t)\le t^2/n\}$, therefore
the last expression tends to zero, uniformly in $a>0$, as $T$ tends to infinity.
This, combined with the classical functional limit theorem corresponding to \eqref{ass:CLT}
and the fact that for a fixed $T>0$, supremum on $[0,T]$ is a continuous map, yields the thesis of \autoref{thm:Seva}.

On the other hand, as a consequence of the Wiener--Hopf factorization (see \cite[Chapter 6]{Kyprianou06}),
with $X^{(a)}_t=X_t-at$, the LST of $\bx a$ is given by
\[
\ee e^{-\lambda\bx a}=
\exp\left(-\int_0^\infty\frac{1}{t}\,\ee\left(1-e^{-\lambda(X_{n(a)t}-an(a) t)},\, X_{n(a)t}-an(a) t>0\right) \,dt \right).
\]
Plugging in $\lambda=\Delta(a)s$ for $s>0$, from \eqref{ass:CLT} and \eqref{eq:defna} it follows that, as $a\downarrow0$, this expression tends to
\[
\ee e^{-s\lr}=
\exp\left(-\int_0^\infty\frac{1}{t}\, \ee\left(1-e^{-s (\LA_t-t)},\, \LA_t-t>0\right)\,dt \right),
\]
the LST of $\lr=\sup_{t\ge0}(\LA_t-t)$,
provided that we can interchange the limit with the integral. 
This follows by using the dominated convergence theorem.
For big values of $t$, say $t>T$ and some $C_3,C_4>0$, we can estimate the integrand by 
(cf. \eqref{eq:Pruitt} and \eqref{eq:estimate})
\[
\frac{1}{t}\pp\left( X_{n(a) t} >a n(a) t\right)\le C\frac{V(an(a)t)}{a^2n(a)t^2}
\le C_3 t^{-\alpha} \frac{V(d(n(a)))}{c^2(n(a))}n(a)\le C_4 t^{-\alpha}.
\]
For $t\le T$ and some $C_5>0$, one can simply bound the integrand by (cf. \eqref{eq:Pruitt})
\[
C_5s t^{1/\alpha-1}\ee(\LA_1,\LA_1>0).
\]
\end{remark}

\subsection{Heavy-traffic invariance principle}
\label{subsec:HTIP}
A general principle called {\it heavy-traffic invariance principle}
has been established by \citet{Szczotka03}; see also \cite{Czystolowski07,Czystolowski10,Szczotka04}.
This principle asserts
under what condition one can infer the limiting distributions of maxima of random walks
from functional limit theorems. According to \autoref{thm:main} this principle
can be also reformulated to the L\'evy setting. Therefore we conclude the paper 
with the following theorem:
\begin{proposition}[Heavy-traffic invariance principle]
For a family of L\'evy processes $\{X^{(a)}_t:t\ge 0\}$ denote
$\mu^{(a)}=\ee X^{(a)}_1<0$ and assume that $\mu^{(a)}\uparrow 0$. Moreover,
assume that there exist functions $d(\cdot)$
and $\Delta(\cdot)$, such that the following conditions hold:
\begin{enumerate}
\item[(I)] $d(a)\Delta(a)|\mu^{(a)}|\to\beta\in(0,\infty)$;
\item[(II)] $\Delta(a)\{X^{(a)}_{d(a)t}-td(a)\mu^{(a)}:t\ge0\}\dto \{X_t:t\ge0\}$ in $D[0,\infty)$,
where $X$ is a L\'evy process;
\item[(III)] the family $\{\Delta(a)\bx a:a>0\}$ is tight.
\end{enumerate}
Then,
\[
\Delta(a)\sup_{t\ge 0}X^{(a)}_t\dto\sup_{t\ge 0}\left ( X_t-\beta t\right).
\]
\end{proposition}
See \cite[Theorem 2]{Szczotka03} for an extension to processes $X^{(a)}$ with stationary 
increments in the case $X$ is stochastically continuous.
\small
\bibliography{Heavy.Traffic}
\end{document}